\documentstyle[12pt,amssymb]{article}
\def\R{{\Bbb R}}
\def\Z{{\Bbb Z}}
\def\N{{\Bbb N}}
\def\eas{\begin{eqnarray*}}
\def\eeas{\end{eqnarray*}}
\def\lll{\lefteqn}
\def\nn{\nonumber}
\def\eq#1{\begin{equation}\label{#1}}
\def\eeq{\end{equation}}
\def\ea#1{\begin{eqnarray}\label{#1}}
\def\eea{\end{eqnarray}}
\def\la{\label}
\def\re#1{(\ref{#1})}

\def\tr{\mathop{\mbox{\rm tr}}}
\def\a{{\alpha}}
\def\b{{\beta}}
\def\g{{\gamma}}
\def\d{\partial}
\def\th{{\theta}}

\def\l{{\lambda}}

\def\del{{\delta}}
\def\eps{{\varepsilon}}

\def\qed{\hfill$\Box$\\\strut}

\def\<{\langle}
\def\>{\rangle}
\def\div{\mbox{\rm div}}
\def\AA{{\mathbf{A}}}
\def\BB{{\mathbf{B}}}
\def\aa{{\mathbf{a}}}
\def\bb{{\mathbf{b}}}
\def\aas{{\mathbf{\tilde{a}}}}
\def\bbs{{\mathbf{\tilde{b}}}}
\def\As{\widetilde{A}}
\def\sysa#1{$\mbox{\re{sys1}}_{#1}$}
\def\sysb#1{$\mbox{\re{sys2}}_{#1}$}
\def\mm{{\mathbf m}}
\def\bib#1{\bibitem[#1]{#1}}
\def\br#1#2{{\frac{#1}{#2}}}

\newtheorem{theorem}{Theorem}[section]
\newtheorem{lemma}[theorem]{Lemma}

\begin{document}
\sloppy
\title{Conservation laws for even order systems of polyharmonic map type}
\author{Fr\'ed\'eric Louis de Longueville \and Andreas Gastel}
\maketitle

{\abstract \noindent
  Following Rivi\`ere's study of conservation laws for second
  order quasilinear systems with critical nonlinearty and Lamm/Rivi\`ere's
  generalization to fourth order, we consider similar systems of order
  $2m$. Typical examples are $m$-polyharmonic maps. Under natural conditions,
  we find a conservation law for weak solutions on $2m$-dimensional domains.
  This implies continuity of weak solutions.\\[6mm]
}

\section{Introduction}

The regularity of harmonic maps and related systems has been an active
topic of research for decades. The harmonic map system in two dimensions has
a critical nonlinearity, and in general such systems can have solutions
that fail to be continuous. But not so for harmonic maps of surfaces. For
a long time, it was not clear which are the general structural assumptions
for critical nonlinear systems in two dimensions that allow proving
full regularity of solutions.

The question was finally solved by Rivi\`ere
\cite{Ri} in 2007. The harmonic map system for harmonic maps
$u:\R^2\supset U\to N$ is usually written as $\Delta u-\tr A(u)(du,du)=0$,
where $A$ is the second fundamental form of the target manifold
$N\subset\R^n$. Rivi\`ere re-interpreted the system as a {\em linear\/}
system
\[
  -\Delta u=\<\Omega,du\>,
\]
where $\Omega$ is a matrix-valued
$1$-form depending on $x$, $u$, and $du$. The point is that
we can ``forget'' the
dependence of $\Omega$ on $u$ and $du$ once we have identified the important
structural assumptions we may impose on $\Omega$. The price one has to
pay is that it has relatively little regularity. Since it depends linearly
on $Du$, it will be a priori only in $L^2$. The tiny bit of extra structure
for $\Omega$ from the harmonic map system is that it can be decomposed
according to $\Omega=d\eta+F$, where $F\in W^{1,2,1}$ (a Sobolev-Lorentz
space slightly smaller than $W^{1,2}$) and $\eta\in W^{1,2}$, but with
values in the skew-symmetric matrices. This antisymmetry in the least
regular term, together with the slight improvement from $W^{1,2}$ to
$W^{1,2,1}$ for the remaining term, are what is needed to prove continuity
of weak solutions. We do not need to know the system is about harmonic maps ---
it is just the structure of $\Omega$ that enters the proof.

The continuity then follows from the fact that the system can be rewritten
as a conservation law
\[
  d({*}A\,du-({*}B)\wedge du)=0,
\]
provided that one can find functions $A\in W^{1,2}\cap L^\infty(U,GL(n))$,
$B\in W^{1,2}(U,so(n)\otimes\wedge^2\R^2)$ satisfying
\[
  dA-A\Omega=-\del B.
\]

For fourth order systems, the situation is quite similar. Regularity for
weak solutions of critically nonlinear systems in four dimensions
cannot be expected to hold in general. There are several variants of
{\em biharmonic\/} map equations, for which regularity is known to hold.
And in 2008, Lamm and Rivi\`ere \cite{LR}
identified structural aspects of the different
biharmonic systems that are responsible for continuity of weak solutions.
Again, they rewrite the systems in the form of a linear equation
\[
  \Delta^2u=\Delta\<V,du\>+\del(w\,du)+\<W,du\>
\]
where the coefficient functions $V\in W^{1,2}$, $w\in L^2$ and $W\in W^{-1,2}$
a priori have little regularity. But again, the least regular coefficient
$W$ is actually slightly better, namely $W=d\eta+F$ with  $F\in L^{4/3,1}$
and $\eta\in L^2$ and $\eta$ skew-symmetric. Note that we have
$L^{4/3}\hookrightarrow W^{-1,2}$ and that $F$ is slightly better.

Again, Lamm and Rivi\`ere find a conservation law
\[
  \del[d(A\Delta u)-2dA\,\Delta u+\Delta A\,du-Aw\,du+dA\<V,du\>
  -Ad\<V,du\>-\<B,du\>]=0,
\]
provided that there are functions $A\in W^{2,2}\cap L^\infty(U,GL(n))$,
$B\in W^{1,\frac43}(U,\R^{n\times n}\otimes\wedge^2\R^4)$ for which
\[
  \Delta dA+(\Delta A)V-(dA)w+AW=\del B.
\]

In both cases, the hard part of the regularity proof is to show that the
auxiliary functions $A$ and $B$ exist under a smallness condition that
can be satisfied by working on sufficiently small balls.

In this paper, we generalize the results from \cite{Ri} ($m=1$) and
\cite{LR} ($m=2$) to systems of order $2m$ on a $2m$-dimensional domain.
The role that harmonic maps and biharmonic maps played so far, will now
be taken by (extrinsically or inrinsically) {\em polyharmonic maps. }They
are critical points
$u:\R^{2m}\supset U\to N\subset\R^n$ of the {\em extrinsic\/} polyharmonic
energy $E_m(u):=\int_U|D^mu|^2\,dx$ or the {\em intrinsic\/} polyharmonic
energy $E_m^i(u):=\int_U|\nabla^{m-1}du|^2\,dx$. Again, weak solutions of
the corresponding Euler-Lagrange equations are smooth, while the same
cannot be expected for solutions of a general semilinear equation of
order $2m$ with critical nonlinearity. The structure needed to prove continuity
of weak solutions along the lines of \cite{Ri} and \cite{LR} turns out
to be quite similar to what we learned above. This time, the formally
linear equation reads
\[
  \Delta^m u=\sum_{k=0}^{m-1}\Delta^k\<V_k,du\>+\sum_{k=0}^{m-2}\Delta^k\del(w_k\,du),
\]
where again we have coefficient functions that in the model case depend
on $u$ and its derivatives, they are
\eas
  w_k&\in&W^{2k+2-m,2}(B^{2m},\R^{n\times n})\qquad
    \mbox{ for }k\in\{0,\ldots,m-2\},\\
  V_k&\in&W^{2k+1-m,2}(B^{2m},\R^{n\times n}\otimes\wedge^1\R^{2m})\qquad
    \mbox{ for }k\in\{0,\ldots,m-1\},\mbox{ where }\\
  V_0&=&d\eta+F,\quad \eta\in W^{2-m,2}(B^{2m},so(n)),
    \quad F\in W^{2-m,\br{2m}{m+1},1}(B^{2m},\R^{n\times n}\otimes\wedge^1\R^{2m}).
\eeas
Note that the regularity assumption on the least regular coefficient
$V_0$ is quite analogous to those made above. It is slightly better than the
$V_0\in W^{2-m,2}$ one would have in the general critical nonlinearity,
in the sense that it decomposes in a skew-symmetric term and a term enjoying
a tiny bit of extra (Lorentz-)regularity. What we can gain from this are
again a conservation law
(see Theorem \ref{th} below for the precise statement) and continuity of
weak solutions (see Theorem \ref{co}). 

However, the setting is slightly less pleasant here than in the
previous cases, because we no longer have reduced everything to coefficients
that are functions. Now $\eta$ and several other coefficients are
in Sobolev spaces of negative order, which means they only make sense
as distributions. This results in a conservation law that depends on
a function $A\in W^{m,2}\cap L^\infty$ and a distribution $B\in W^{2-m,2}$
that satisfy some auxiliary equation. It is this difference between
functions and distributions that makes our considerations more technical
than those of \cite{LR}. Rather than splitting the auxiliary equation
for $A$ and $B$ into a system of two equations, we have to solve a much
larger system, roughly to determine all $B_\a$ in a decomposition
$B=\sum_{|\a|\le m-2}D^\a B_\a$.

Apart from that extra complication from dealing with negative
differentiability, our methods are reasonably close to Rivi\`ere's
and Lamm/Rivi\`ere's. In particular, we use Uhĺenbeck's gauge theorem, Lorentz
spaces and Hodge decomposition.

Regularity of weakly polyharmonic maps has already been proven by Scheven and
the second author in \cite{GS}. It is not our primary focus here to re-prove
that theorem, but to study the problem under minimal structural conditions
on the equations as in the works by Rivi\`ere and Lamm/Rivi\`ere.

This paper is based on the first author's Ph.D. thesis. The work was
supported by DFG grant GA1428/5-1.

\section{Preliminaries}
\subsection{Facts about Lorentz spaces}\la{sec:ls}

Let $U\subset\R^d$ be a bounded domain with smooth boundary.
Lorentz spaces $L^{p,q}(U)$ with $p,q\in[1,\infty]$ interpolate between the
Lebesgue spaces $L^p(U)$ in such a way that $L^{p,p}(U)=L^p(U)$ and
$L^{p,1}(U)\subset\ldots\subset L^{p,p}(U)\subset\ldots\subset L^{p,\infty}(U)$
for all $p$. We refer to Ziemer's book \cite[\S\,1.8 and \S\,2.10]{Zi}
for definitions and basic properties. We are going to use embedding properties
and Calderon-Zygmund estimates for Lorentz spaces rather than the actual
definition.

For $k\in\N$, we denote by $W^{k,p,q}(U)$ the Lorentz-Sobolev space of all
functions $f$ on $U$ which are weakly differentiable $k$ times with
$D^\a f\in L^{p,q}(U)$ for all $|\a|\le k$. The $W^{k,p,q}$-norm is the
sum of the $L^{p,q}$-norms of the $D^\a f$.

For $k\in\Z_{<0}$, $W^{k,p,q}$-spaces are often defined as duals of
Sobolev-Lorentz spaces. But doing so, we would miss the spaces for
$q=1$, which we will use throughout the paper. We therefore use the alternative
definition, which is equivalent to the dual one in the $q>1$ and $p>1$
cases. For $k\in\Z_{<0}$, $p,q\in[1,\infty]$, we define $W^{k,p,q}(U)$
to be the space of all distributions on $U$ of the form
$f=\sum_{|\a|\le-k}D^\a f_\a$ with $f_\a\in L^{p,q}(U)$. The corresponding
norm is defined as $\|f\|_{W^{k,p,q}(U)}:=\inf\sum_{|\a|\le-k}\|f_\a\|_{L^{p,q}(U)}$,
where the infimum is taken over all decompositions of $f$ as given in the
definition.

Here are some facts about (Sobolev-)Lorentz spaces we are going to use.

{\bf (1) Sobolev-Lorentz embeddings. }They have been proven by Tartar \cite{Ta}
for Sobolev-Lorentz spaces on $\R^d$.
Using extension theorems which hold analogous
to the ones for usual Sobolev spaces, we easily carry them over to smooth
bounded domains. If $k,\ell\in\Z$, $1<p<\frac{d}{\ell}$ and $1\le q\le\infty$,
then every $f\in W^{k,p,q}(U)$ is also in $W^{k-\ell,\br{dp}{d-\ell p},q}(U)$ with
\[
  \|f\|_{W^{k-\ell,\br{dp}{d-\ell p},q}(U)}\le C\|f\|_{W^{k,p,q}(U)}\,.
\]

{\bf (2) Multiplication theorems. }Assume $f\in W^{p,q}(U)$ and
$g\in W^{p',q'}(U)$, where $1<p,p'<\infty$ with $\frac1{p}+\frac1{p'}<1$ and
$1\le q,q'\le\infty$. Then
\[
  \|fg\|_{L^{r,s}}\le C\|f\|_{L^{p,q}}\|g\|_{L^{p',q'}}
\]
if $\frac1r=\frac1{p}+\frac1{p'}$ and $\frac1s\le\frac1{q}+\frac1{q'}$,
$s\ge1$. In the boundary case $\frac1{p}+\frac1{p'}=1$, we have
\[
  \|fg\|_{L^1}\le\|f\|_{L^{p,q}}\|g\|_{L^{p',q'}}
\]
whenever $\frac1{q}+\frac1{q'}\le 1$. See \cite[Thm.\ 3.4, 3.5]{ON}
for proofs. From these theorems together with the Sobolev-Lorentz embeddings
above, we also have Sobolev-Lorentz multiplication theorems. Suppose
$f\in W^{k,p,q}(U)$ and $g\in W^{z,p',q'}$ with $k\in\N\cup\{0\}$,
$z\in\Z$, $1<p,p'<\infty$, $1\le q,q'<\infty$. If $\frac1{p}+\frac1{p'}\le1$,
$k\ge|z|$, % (WIRKLICH?)
and $kp<d$, then $fg\in W^{z,x,y}(U)$ for
$x:=\frac{dpp'}{d(p+p')-kpp'}$ and $\frac1y:=\min\{1,\frac1q+\frac1{q'}\}$,
and
\eq{zxy}
  \|fg\|_{W^{z,y,x}}\le C\|f\|_{W^{k,p,q}}\|g\|_{W^{z,p',q'}}\,.
\eeq
The assertion continues to hold in the case $kp=d$ (and then with $x=p'$
and $y=q'$) if we additionally assume $f\in L^\infty$.

{\bf (3) $L^p$ theory.} The same way we have $W^{k+2,p}$ estimates for
solutions of $\Delta u=f\in W^{k,p}$ with suitable boundary data, we also
have $W^{k+2,p,q}$-estimates if $\Delta u=f\in W^{k,p,q}$. This is proven
by standard interpolation arguments and was used extensively in
Rivi\`ere's and Lamm's papers. Similarly, Hodge decomposition theorems
in $W^{k,p}$ extend to $W^{k,p,q}$. Using our definition of $W^{k,p,q}$ for
$k<0$, we can even extend most of the results to those spaces by
applying them to a decomposition $f=\sum_{|\a|\le-k}D^\a f_\a$. We have to
be careful with boundary values, though, imposing them on the ``partial
solutions'' $u_\a$ of $\Delta u_\a=f_\a$ and finally defining
$u:=\sum_{|\a|\le-k}D^\a u_\a$. Reducing any solution of $\Delta f=g$ to
a solution of such boundary value problems using cutoff functions,
it is standard to prove the following.

\begin{lemma}[some basic $L^{p,q}$-theory]\la{Lpq}
Assume $k\in\Z$, $p\in(1,\infty)$, $q\in[1,\infty)$.
Given $g\in W^{k,p,q}(B_1)$, any weak solution $f$ of $\Delta f=g$
is in $W^{k+2,p,q}(B_{1/2})$, and we have
\[
  \|f\|_{W^{k+2,p,q}(B_{1/2})}\le c\|g\|_{W^{k,p,q}(B_1)}
\]
with a constant $c$ depending only on $k,p,q$ and the domain dimension.
\end{lemma}

{\bf (4) A nontrivial Hodge-style Lemma. }The Lemma that follows is taken
from \cite[Lemma A.2]{Ri} for $m=1$ and \cite[Lemma A.2]{LR} for $m=2$. Their
proof carries over immediately. Even though the assertion looks little
surprising in view of Hodge decompositions, it must be emphasized that
this is one of the least trivial ingredients in this paper. The proof
does involve the Coifman-Lions-Meyer-Semmes-Lemma
\cite{CLMS}, Hardy-BMO-duality, and the fact that
the Neumann problem for the Poisson equation with right-hand side in
${\mathcal H}^1(\R^d)$ can be solved in $W^{1,\br{d}{d-1},1}(\R^d)$. That
is, the Lemma uses all the facts about ``exotic'' function spaces
that have been found to be important for harmonic map theory
by H\'elein \cite{He} and many others.

\begin{lemma}\la{TODO}
  Assume $m\in\N$. There is $\eps>0$ such that for all
  $Q\in W^{m,2}\cap L^\infty(B^{2m},SO(n))$ satisfying $Q-I\in W^{m,2}_0$
  and $\|dQ\|_{L^{2m}(B^{2m})}<\eps$, the boundary value problem
  \[\left\{\begin{array}{ll}
    d[(\del C)Q]&\;=\;0\qquad\mbox{ on }B^{2m},\\
    dC&\;=\;0\qquad\mbox{ on }B^{2m},\\  
    C_N&\;=\;0\qquad\mbox{ on }\d B^{2m}
  \end{array}\right.\]
  in $W^{1,\br{2m}{2m-1}}$ has the trivial solution $C\equiv0$ as the only
  solution.
\end{lemma}
  
{\bf (5) Continuity of $W^{k,d/k,1}$-functions. }An important point in proving
continuity of our weak solutions will be to prove that they are in
$W^{m+1,\frac{2m}{m+1},1}$ and then use that this space embeds into
$C^0$. This comes from an embedding which generalizes the well-known
embedding $W^{d,1}\to C^0$. It was used, but not explicitly stated,
in \cite[Section 2.3]{LR} for $d=4$, $k=2$. Our proof uses the arguments
presented there.

\begin{theorem}\la{EmbC0}
  Let $\Omega\subset\R^d$ be any bounded domain of class $C^{k+1}$, and assume
  $k\in\{1,\ldots,d\}$. Then any map $u\in W^{k,d/k,1}(\Omega)$ is also in
  $C^0(\Omega)$ and obeys the estimate
  \[
    \|u\|_{L^\infty(\Omega)}\le c\|u\|_{W^{k,d/k,1}(\Omega)}\,.
  \]
\end{theorem}

{\bf Proof.}
Since the embedding $W^{d,1,1}=W^{d,1}\hookrightarrow C^0$ is well-known, we may
assume $k<d$. Assume that $k$ is odd. Let $G$ be the fundamental solution of
$\Delta^{(k+1)/2}$ on $\R^d$, $f:=dG$, which means $f(x)=c|x|^{k-d-1}x$ and
hence $f\in L^{d/(d-k),\infty}(\Omega)=L^{d/k,1}(\Omega)^*$. Extend $u$ to a
function $u\in W^{k,d/k,1}(\R^d)$ with norm controlled by the one on the original
$u$. For any $x\in\R^d$, we have
\eas
  |u(x)|&=&\Big|\int_{\R^d}\Delta^{(k+1)/2}G(x-y)u(y)\,dx\Big|
  \;=\;\Big|\int_{\R^d}\<f(x-y),d\Delta^{(k-1)/2}u(y)\>\,dy\Big|\\
  &\le&\|f\|_{L^{d/(d-k),\infty}(\R^d)}\,\|d\Delta^{(k-1)/2}u\|_{L^{d/k,1}(\R^d)}
  \;\le\;c\|u\|_{W^{k,d/k,1}(\Omega)},
\eeas
which proves the $L^\infty$-estimate. Now approximate $u$ in $W^{k,d/k,1}(\R^d)$
by functions $u_j\in C^\infty_c(\R^d)$. By the estimate just proven,
$\|u_j-u_i\|_{L^\infty(\R^d)}\le c\|u_j-u_i\|_{W^{k,d/k,1}(\R^n)}$, which shows that
$(u_j)$ is a uniform Cauchy sequence, implying continuity of $u$.

If $k$ is even, a similar reasoning uses the fundamental solution $G$ for
$\Delta^{k/2}\del$ and $f:=\del G$.\qed

\subsection{Uhlenbeck decomposition}

We will need a suitable adaptation of Uhlenbeck's gauge theorem
\cite[Theorem 1.3]{Uh}
in the spirit of Rivi\`ere's reinterpretation \cite[Lemma A.4]{Ri}, cf.\ 
\cite[Theorem A.5]{LR} for the $m=2$ case. See also \cite{GZ} for a discussion
of its analytic aspects.

\begin{theorem}\la{uhl}
Assume that $m,n\in\N$ and $B_r\subset\R^{2m}$ is a ball of radius $r$.
Then there is $\eps>0$ such that for all 
$\Omega\in W^{m-1,2}(B_r,so(n)\otimes\wedge^1\R^{2m})$ satisfying
$\|\Omega\|_{W^{m-1,2}(B_r)}<\eps$, there are functions 
$P\in W^{m,2}(B_{r/2},SO(n))$ and
$\Phi\in W^{m,2}(B_{r/2},so(n)\otimes\wedge^2\R^{2m})$ such that
\[
  \Omega=P\,dP^{-1}+P\,\del\Phi\,P^{-1}
\]
holds on $B_{r/2}$. Moreover, we have the estimate
\[
  \|dP\|_{W^{m-1,2}(B_{r/2})}+\|\del\Phi\|_{W^{m-1,2}(B_{r/2})}\le c\|\Omega\|_{W^{m-1,2}(B_r)}.
\]
\end{theorem}

{\bf Proof.} The proof is given in the references above, except for the
estimates which we find there only for $\|dP\|_{W^{1,2}}$ and
$\|\del\Phi\|_{W^{1,2}}$. We can have $d\Phi=0$ there, and bootstrapping
the estimates through the equations $d\Phi=0$ and
$P^{-1}dP^{-1}+P^{-1}\Omega P=\del\Phi$ easily gives the higher order
estimates.\qed

\section{Polyharmonic maps and their equations}

\subsection{Extrinsically polyharmonic maps}

Extrinsically polyharmonic maps $u:\R^n\supset U\to N$
are critical points of $\int|D^mu|^2\,dx$.
The most natural way to write down the Euler-Lagrange equation for that
functional is
\[
  \Delta^mu\perp N.
\]
To make this an explicit system of differential equations, we use local
smooth orthonormal bases $\{\nu_i(u(x))\}$ of $T_{u(x)}^\bot N$ and can
therefore write
\[
  \Delta^mu=\sum_i\l_i(\nu_i\circ u),
\]
with Lagrange multipliers $\l_i$ that we are going to determine. 
Multiplying with $(\nu\circ u)$ and using Leibniz' rule iteratively, we find
\eas
  \l_i&=&\<\Delta^mu,\nu_i\circ u\>\\
  &=&\Delta^m\<u,\nu_i\circ u\>-\sum_{t=0}^{m-1}\Delta^t\<d\Delta^{m-t-1}u,J_idu\> 
  -\sum_{t=0}^{m-1}\Delta^t\del\<\Delta^{m-t-1}u,J_idu\>,
\eeas
where here $J_i:=d\nu_i\circ u$. Expanding the innermost $\Delta=\del d$,
the first term on the right-hand side becomes
\[
  \Delta^m\<u,\nu_i\circ u\>=\Delta^{m-1}\del[\<du,\nu_i\circ u\>+\<u,J_idu\>]
  =\Delta^{m-1}\del\<u,J_idu\>,
\]  
since $du$ is tangential to $N$ and $\nu_i\circ u$ is normal. Hence
\eq{eq1}
  \l_i\nu_i=-\sum_{t=0}^{m-1}\Delta^t\<d\Delta^{m-t-1}u,J_idu\>\nu_i 
  -\sum_{t=0}^{m-2}\Delta^t\del\<\Delta^{m-t-1}u,J_idu\>\nu_i,
\eeq
where we have abbreviateted $\nu_i\circ u$ by $\nu_i$ and will keep doing
so. The $t=0$ term in the first sum can be rewritten as
\eq{wterm}
  -\<\Delta^{m-1}du,J_idu\>\nu_i=-\<\nu_i\otimes J_i^T\Delta^{m-1}du,du\>
\eeq
while for $t>0$ a long but straightforward calculation shows
\eas
  \lll{-\Delta^t\<\Delta^{m-t-1}du,J_idu\>\nu_i\;=\;-\Delta^t\<\nu_i\otimes J_i^T
    \Delta^{m-t-1}du,du\>}\\
  &&{}+\sum_{j=1}^t\Delta^{t-j}\Big\<J_i^T\Delta^{j-1}d\<\Delta^{m-t-1}du,J_idu\>,
    du\Big\>+\sum_{j=1}^t\Delta^{t-1}\del[\Delta^{j-1}\<\Delta^{m-t-1}du,
    J_idu\>J_idu].
\eeas
The $t=0$ term in the second sum of \re{eq1} becomes
\[
  -\del\<\Delta^{m-1}u,J_idu\>\nu_i=-\del\<\nu_i\otimes J_i^T\Delta^{m-1}u
    ,du\>+\Big\<J_i^T\<\Delta^{m-1}u,J_idu\>,du\Big\>,
\]
while for $t>0$ we have 
\eas
  \lll{-\Delta^t\del\<\Delta^{m-t-1}u,J_idu\>\nu_i\;=\;-\Delta^t\del\<\nu_i
    \otimes J_i^T\Delta^{m-t-1}u,du\>+\Delta^t\Big\<\<J_i\Delta^{m-t-1}u,
    J_idu\>,du\Big\>}\\
  &&{}+\sum_{j=1}^t\Delta^{t-j}\Big\<J_i^T\Delta^{j-1}d\del\<\Delta^{m-t-1}u,
    J_idu\>,du\Big\>+\sum_{j=1}^t\Delta^{t-j}\del\Big\<\Delta^{j-1}\del\<
    \Delta^{m-t-1}u,J_idu\>,J_idu\Big\>.
\eeas
Inserting that into \re{eq1}, we find
\eq{eq2}
  \Delta^m u=\sum_{k=0}^{m-1}\Delta^k\<V_k,du\>+\sum_{k=0}^{m-2}\Delta^k\del(w_k\,du)
\eeq
with (abbreviating $d^\ell:=d\del d...$ and $\del^\ell:=\del d\del...$ with
$\ell$ letters each)
\ea{eqV}
  V_k&:=&\sum_i\Big(-(\nu_i\circ u)\otimes J_i^T\Delta^{m-k-1}du
    +\sum_{t=k+1}^{m-1}J_i^T\Delta^{t-k-1}d\<\Delta^{m-t-1}du,J_idu\>\nn\\
  &&\qquad\quad+\sum_{t=k}^{m-2}J_i^Td^{2t-2k}\<\Delta^{m-t-1}u,J_idu\>\Big)\\
\la{eqw}
  w_k&:=&\sum_i\Big(-(\nu_i\circ u)\otimes J_i^T\Delta^{m-k-1}u
    +\sum_{t=k+1}^{m-1}J_i^T\Delta^{t-k-1}\<\Delta^{m-t-1}du,J_idu\>\nn\\
  &&\qquad\quad+\sum_{t=k}^{m-2}J_i^T\del^{2t-2k}\<\Delta^{m-t-1}u,J_idu\>\Big)
\eea
For weakly polyharmonic mappings, we must assume $u\in W^{m,2}$, and then
some of the terms in $V_k$ and $w_k$ exist only in the sense of distributions.
But it is easily checked using Section \ref{sec:ls}\,(1)--(2)
that the terms do not involve products of
distributions (which would be undefined), and that
\eas
  V_k\in W^{2k+1-m,2}&&\mbox{ for }k=0,\ldots,m-1,\\
  w_k\in W^{2k+1-m,2}&&\mbox{ for }k=0,\ldots,m-2.
\eeas
But we will need that $V_0$ is better than that, and it actually can be
improved. This
is the important extra structure that Rivi\`ere discovered to be essential
in \cite{Ri} for second-order systems (and in \cite{LR} for fourth order).
We can modify \re{wterm} by adding a term that is $0$ because of
$\nu_i\perp du$,
\eas
  \lll{-\sum_i\<\Delta^{m-1}du,J_idu\>\nu_i\;=\;
    \sum_i\<J_i^T\Delta^{m-1}du\otimes\nu_i-\nu_i\otimes J_i^T\Delta^{m-1}du,du\>}
  \nn\\
  &=&\sum_i\<d[J_i^T\Delta^{m-1}u\otimes\nu_i-\nu_i\otimes J_i^T\Delta^{m-1}u],du\>
  +\<R,du\>\nn\\
  &=:&\<d\eta,du\>+\<R,du\>,
\eeas
where here $\eta$ takes its values in the set $so(n)$ of skew-symmetric
$n\times n$-matrices, and $\eta\in W^{2-m,2}$. The remaining terms in $V_0$,
including new ones from $R$, are all products of at least two
$W^{...,2}$-functions, hence they are in some Sobolev-Lorentz space
$W^{...,...,1}$. A case-by-case inspection of all summands shows that our modified
$V_0$ splits into
\eq{regV0}
  V_0=d\eta+F,\quad\eta\in W^{2-m,2}(U,so(n)),\quad
  F\in W^{2-m,\br{2m}{m+1},1}(U,\R^{n\times n}\otimes\wedge^1\R^{2m}).
\eeq  
This is as much structure as we will need for our regularity theory.

\subsection{Intrinsically polyharmonic maps}

Intrinsically polyharmonic maps are critical points of
$\int|\nabla^{m-1}du|^2\,dx$. Since
\[
  \nabla^{m-1}du-D^mu=\sum_{\mm\in M_1}(X_\mm\circ u)(D^{m_1}u,\ldots,D^{m_\ell}u)
\]
for suitable multilinear forms $X_\mm$ and
\[
  M_1:=\Big\{\mm=(m_1,\ldots,m_\ell)\in\bigcup_{\ell=2}^m\N^\ell:m_1+\ldots
  +m_\ell=m,\,m_1\ge\ldots\ge m_\ell\Big\},
\]
the difference between the extrinsic and the intrinsic integrand can be written
as
\[
  |\nabla^{m-1}du|^2-|D^mu|^2=\sum_{\mm\in M_2}(Y_\mm\circ u)
  (D^{m_1}u,\ldots,D^{m_\ell}u),
\]
with new multilinear forms $Y_\mm$ and
\[
  M_2:=\Big\{\mm=(m_1,\ldots,m_\ell)\in\bigcup_{\ell=3}^m\{1,\ldots,m\}^\ell
  :m_1+\ldots+m_\ell=2m,\, m_1\ge\ldots\ge m_\ell\Big\}.
\]
This means that the Euler-Lagrange equation for intrinsically polyharmonic maps
reads (somewhat symbolically, since we do not indicate precisely which
contractions are performed by the div-operators)
\[
  \Delta^mu+\sum_{\mm\in M_3}\div^{m_1}((Z_\mm\circ u)(D^{m_2}u,\ldots,D^{m_\ell}u)
  \perp N,
\]
with yet more multilinear forms $Z_\mm$ and
\[
  M_3:=\Big\{\mm=(m_1,\ldots,m_\ell)\in\bigcup_{\ell=3}^m\{1,\ldots,m\}^\ell
  :m_1+\ldots+m_\ell=2m,\, m_2\ge\ldots\ge m_\ell\Big\}.
\]
It is now straightforward to check that the additional terms contribute to
more summands for $V_0,\ldots,V_{m-1},w_0,\ldots,w_{m-2}$ (in the same spaces)
in the previous
calculation, and that $\ell\ge3$ in each summand implies that no derivatives
of order $2m-1$ will be involved in the terms contributing to $V_0$. Hence
none of those will contribute to $d\eta$. This means that the Euler-Lagrange
equation for intrinsically polyharmonic maps can be written in the form
\re{eq2} with exactly the same regularity conditions for the coefficient
functions as in the extrinsic case.

\subsection{Further equations}

The equations allowed here are by no means restricted to a geometric or
variational context. Equations of the form \re{eq2} do not refer to
a target manifold $N$ explicitly. For example, they include equations like
the ``fake polyharmonic equation''
\[
  \Delta^mu+|Du|^{2m}u=0,
\]
whose coefficient functions fulfill the required estimates if $u$ is
bounded.

\section{A conservation law}

The following theorem reformulates a rather general system of order
$2m$ as a conservation law, in the spirit of Rivi\`ere's original idea
\cite[Theorems I.3 and I.4]{Ri} ($m=1$) and Lamm/Rivi\`ere's fourth order
generalization \cite[Theorems 1.3 and 1.5]{LR} ($m=2$).
We work on the unit ball $B^{2m}$ of $\R^{2m}$ and later
use a scaling argument.

\begin{theorem}\la{th}
Assume $m\ge3$, $n\in\N$. Let coefficient functions be
given as
\eas
  w_k&\in&W^{2k+2-m,2}(B^{2m},\R^{n\times n})\qquad
    \mbox{ for }k\in\{0,\ldots,m-2\},\\
  V_k&\in&W^{2k+1-m,2}(B^{2m},\R^{n\times n}\otimes\wedge^1\R^{2m})\qquad
    \mbox{ for }k\in\{0,\ldots,m-1\},\mbox{ where }\\
  V_0&=&d\eta+F,\quad \eta\in W^{2-m,2}(B^{2m},so(n)),
    \quad F\in W^{2-m,\br{2m}{m+1},1}(B^{2m},\R^{n\times n}\otimes\wedge^1\R^{2m}).
\eeas
We consider the equation
\eq{eq}
  \Delta^m u=\sum_{k=0}^{m-1}\Delta^k\<V_k,du\>+\sum_{k=0}^{m-2}\Delta^k\del(w_k\,du).
\eeq
For this equation, the following statements hold.

(i) Let
\eas
  \th&:=&\sum_{k=0}^{m-2}\|w_k\|_{W^{2k+2-m,2}(B^{2m})}
    +\sum_{k=1}^{m-1}\|V_k\|_{W^{2k+1-m,2}(B^{2m})}\\
  &&\quad{}+\|\eta\|_{W^{2-m,2}(B^{2m})}+\|F\|_{W^{2-m,\br{2m}{m+1},1}(B^{2m})}\,.
\eeas
There is $\th_0>0$ such that whenever $\th<\th_0$, there are a function 
$A\in W^{m,2}\cap L^\infty(B_{1/4};GL(n))$ and a distribution
$B\in W^{2-m,2}(B_{1/4},\R^{n\times n}\otimes\wedge^2\R^{2m})$ that solve
\eq{aux}
  \Delta^{m-1}dA+\sum_{k=0}^{m-1}(\Delta^kA)V_k-\sum_{k=0}^{m-2}(\Delta^kdA)w_k
  =\del B.
\eeq
% (ABSCHAETZUNG DURCH $\sqrt{\th}$?)

(ii) A function $u\in W^{m,2}(B_{1/2},\R^n)$ solves \re{eq} weakly on
$B_{1/4}$ if and only if it is a distributional solution of the conservation
law
\ea{cons}
  0&=&\del\Big[\sum_{\ell=0}^{m-1}(\Delta^\ell A)\Delta^{m-\ell-1}du
    -\sum_{\ell=0}^{m-2}(d\Delta^\ell A)\Delta^{m-\ell-1}u\nn\\
  &&\qquad{}-\sum_{k=0}^{m-1}\sum_{\ell=0}^{k-1}(\Delta^\ell A)\Delta^{k-\ell-1}
      d\<V_k,du\>
    +\sum_{k=0}^{m-1}\sum_{\ell=0}^{k-1}(d\Delta^\ell A)\Delta^{k-\ell-1}\<V_k,du\>\nn\\
  &&\qquad{}-\sum_{k=0}^{m-2}\sum_{\ell=0}^{k-2}(\Delta^\ell A)d\Delta^{k-\ell-1}
      \del(w_k\,du)
    +\sum_{k=0}^{m-2}\sum_{\ell=0}^{k-2}(d\Delta^\ell A)\Delta^{k-\ell-1}
      \del(w_k\,du)\nn\\
  &&\qquad{}-\<B,du\>\Big].
\eea
(Here $d\Delta^{-1}\del$ means the identity map.)

(iii) Every weak solution of $\re{eq}$ on $B^{2m}$ is continuous on
$B_{1/16}$ if the smallness condition $\th<\th_0$ holds.
\end{theorem}

{\bf Proof of Theorem \ref{th} (ii).} This can be done by direct calulation.
A line-by-line calculation of the terms in \re{cons} gives
\eas
  \del[...]&=&\<d\Delta^{m-1}A,du\>+A\Delta^mu\\
  &&{}-\sum_{k=0}^{m-1}\{A\Delta^k\<V_k,du\>-(\Delta^kA)\<V_k,du\>\}\\
  &&{}-\sum_{k=0}^{m-2}\{\<d\Delta^kA,w_k\,du\>+A\Delta^k\del(w_k\,du)\}\\
  &&{}-\<\del B,du\>\\
  &=&A\Big\{\Delta^m u-\sum_{k=0}^{m-1}\Delta^k\<V_k,du\>-\sum_{k=0}^{m-2}
    \Delta^k\del(w_k\,du)\Big\}.
\eeas
Now $\re{eq}$ means that the last $\{...\}$ vanishes, and $A$ is invertible
by (i).\qed

{\bf Proof of Theorem \ref{th} (iii).} On $B_{1/4}$, the conservation law
\re{cons} can be rewritten as
\[
  0=\del[A\Delta^{m-1}du+(dA)\Delta^{m-1}u+R]=\Delta(A\Delta^{m-1}u)+\del R,
\]
where here $R\in W^{2-m,\br{2m}{m+1},1}(B_{1/4},\R^n)$. By $L^{p,q}$ theory in
Lorentz-spaces as formulated in Lemma \ref{Lpq},
we first have $A\Delta^{m-1}u\in 
W^{3-m,\br{2m}{m+1},1}$ on $B_{1/8}$. Invertibility of $A$ and $A\in 
W^{m,2}\cap L^\infty$ give $\Delta^{m-1}u\in W^{3-m,\br{2m}{m+1},1}$ by the remark
following \re{zxy}, and therefore
$u\in W^{m+1,\br{2m}{m+1},1}$ on $B_{1/16}$. But the latter space embeds into
$C^0$ by Theorem \ref{EmbC0}, which implies the asserted continuity of $u$.\qed

{\bf Proof of Theorem \ref{th} (i).}
This is the technical part of our paper.
By iteratively solving Neumann problems on $B^{2m}$ (in the case of negative
Sobolev exponents to be understood in the sense described in
Section \ref{sec:ls}\,(3)), we find 
$\Psi\in W^{m,2}(B^{2m},so(n))$ satisfying $\Delta^{m-1}\Psi=-\eta$. Note that
we may assume $\int_{B^{2m}}\eta\,dx=0$ since our $\eta$ enters into the
theorem only in $d\eta$. Letting 
$\Omega:=d\Psi\in W^{m-1,2}(B^{2m},so(n)\otimes\wedge^1\R^{2m})$, we have found 
$\Omega$ satisfying
\[
  \Delta^{m-2}\del\Omega=-\eta,\qquad
  \|\Omega\|_{W^{m-1,2}(B^{2m})}\le c\|\eta\|_{W^{2-m,2}(B^{2m})}\,.
\]
If $\th$ in our theorem has been chosen small enough, then $\Omega$ satisfies
the smallness assumption for Theorem~\ref{uhl}. Then we find corresponding
functions $P$ and $\Phi$ such that
\[
  V_0=d\eta+F=-d\Delta^{m-2}\del\Omega+F=-d\Delta^{m-2}\del(P\,dP^{-1}
  +P\,\del\Phi\,P^{-1})+F
\]
on $B_{1/2}$, where here $dP$ and $\del\Phi$ are small on $B_{1/2}$ in
their respective norms which are bounded by $c\th$. 
We rewrite that further,
\[
  V_0=-Pd\Delta^{m-1}P^{-1}+K,
\]
where here
\eas
  K&=&-d\Delta^{m-2}(\<dP,dP^{-1}\>+\<dP,\del\Phi\>P^{-1}+P\<\del\Phi,dP^{-1}\>)
  +\sum_{j=1}^{2m-3}\<D^jP,D^{2m-3-j}\Delta P^{-1}\>\\
  &\in&W^{2-m,\br{2m}{m+1},1}(B_{1/2}),
\eeas
and the corresponding norm is bounded by $c\th$.
Here and in the sequel, we use $\<D^aU,D^bV\>$ as a symbol for any
bilinear form involving $a$-th derivatives of $U$ and $b$-th derivatives
of $V$, if it can be computed in principle and the explicit form does not
matter.

We use the function $P$ constructed here to substitute $A$ by $(I+\As)P^{-1}$
and find that \re{aux}
transforms to an equation for $(\As,B)$ of the form
\eq{AsB}
  d\Delta^{m-1}\As+\sum_{j=0}^{2m-2}\<D^j\As,K_j\>+K_0=\del B\,P
\eeq
with coefficient functions $K_0,\ldots,K_{m-2}$ bounded by $c\th$
in the norms of the spaces
\[
  K_0\in W^{2-m,\br{2m}{m+1},1}(B_{1/2}),\qquad
  K_j\in W^{j+1-m,2}(B_{1/2})\mbox{ for }j\in\{1,\ldots,m-2\}.
\]
Now assume that instead of $w_k$, $V_k$, $\eta$, and $F$, we have started
with $\psi w_k$, $\psi V_k$, $\psi\eta$, and $\psi F$ for some smooth
$B_{1/4}$-$B_{1/2}$-cutoff function. Then everything is in the same spaces
as before, and continues to be controlled by $c\theta$. The effect on
$P$ and $\Phi$ is $P\equiv I$ and $\Phi\equiv0$ on some neighborhood of
$\d B_{1/2}$. This implies $K_j\equiv0$ there, too, hence we continue to
have an equation of the form \re{AsB}, but with nice boundary behaviour
of $P$ and the $K_j$. Most importantly, $\psi\equiv 1$ on $B_{1/4}$
implies that on $B_{1/4}$ both equations coincide. 

By Lemma \ref{iterier} below, the modified system has a solution on
$B_{1/2}$ that is controlled as needed in the theorem. Hence we have found
a solution $(\As,B)$ to the unmodified equation \re{AsB} on $B_{1/4}$.
Transforming back, we find $(A,B)$ as asserted in part (i) of the Theorem,
with $A$ taking values in $GL(n)$ because of $A=(I+\As)P^{-1}$ and
$\|\As\|_{L^\infty}<c\th^{1/2}$.
Hence the theorem is now proven up to Lemma~\ref{iterier} that follows.\qed

\begin{lemma}\la{iterier}
  Assume we are given functions
  \eas
    P&\in&W^{m,2}(B^{2m},SO(n)),\\
    K_0&\in&W^{2-m,\br{2m}{m+1},1}(B^{2m}),\\
    K_j&\in&W^{j+1-m,2}(B^{2m})\qquad\mbox{ for }j\in\{1,\ldots,2m-2\}
  \eeas
  which are small in their respective spaces,
  \eq{klein}
    \|dP\|_{W^{m-1,2}(B^{2m})}+\|K_0\|_{W^{2-m,\br{2m}{m+1},1}(B^{2m})}
    +\sum_{j=1}^{2m-2}\|K_j\|_{W^{j+1-m,2}(B^{2m})}\le\th
  \eeq
  for some $\th$ to be chosen small enough. Assume that $P\equiv I$ (identity
  matrix) and $K_j\equiv0$ ($0\le j\le 2m-2$) hold on a neighborhood of
  $\d B^{2m}$. Then there exist
  \eas
    A&\in&W^{m,2}\cap L^\infty(B^{2m},\R^{n\times n}),\\
    B&\in&W^{2-m,2}(B^{2m},\R^{n\times n}\otimes\wedge^2\R^{2m})
  \eeas
  such that the equation
  \eq{hilfs}
    d\Delta^{m-1}A+\sum_{j=0}^{2m-2}\<D^jA,K_j\>+K_0=\del B\,P
  \eeq
  holds in the sense of distributions. Moreover,
  \eq{ABklein}
    \|A\|_{W^{m,2}(B^{2m})}+\|A\|_{L^\infty(B^{2m})}+\|B\|_{W^{2-m,2}(B^{2m})}\le c\th^{1/2}.
  \eeq
\end{lemma}

{\bf Proof.}
We are looking for $A$ and $B$ in the form of sums of derivatives
\eas
  A&=&\sum_{|\a|\le m-2}\d^\a A_\a,\qquad A_\a\in W^{2m-1,\br{2m}{2m-1-|\a|},1},\\
  B&=&\sum_{|\a|\le m-2}\d^\a B_\a,\qquad B_\a\in W^{1,\br{2m}{2m-1-|\a|},1}.
\eeas
Note that this even gives $A\in W^{m+1,\br{2m}{m+1},1}$ and
$B\in W^{3-m,\br{2m}{m+1},1}$, but these spaces embed into those stated in
the lemma. By the usual representation of negative Sobolev(-Lorentz)
spaces, we find decompositions of $K_0,K_1,\ldots,K_{m-2}$ of the form
\eas
  K_0&=&\sum_{|\a|\le m-2}\d^\a K_0^\a,\qquad K_0^\a\in L^{\br{2m}{m+1},1},\\
  K_j&=&\sum_{|\a|\le m-1-j}\d^\a K_j^\a,\qquad K_j^\a\in L^2,
\eeas
such that $\sum_\a\|K_j^\a\|\le c\,\|K_j\|$ in their respective norms. For
$1\le j\le m-2$, we then have (with certain integers $c_{\b\g}$)
\ea{l1}
  \<D^jA,K_j\>&=&\sum_{|\a|\le m-2}\,\sum_{|\b|\le m-1-j}\<\d^\a D^jA_\a,
    \d^\b K_j^\b\>\nn\\
  &=&\sum_{\a,\b}\sum_{\g\le\b}\d^\g[c_{\b\g}\<\d^{\b-\g}\d^\a D^jA_\a,K_j^\b\>]\nn\\
  &=:&\sum_{\a,\b}\sum_{\g\le\b}\d^\g\<A,K\>_{j,\a,\b,\g}\,,
\eea
where here $\<A,K\>_{j,\a,\b,\g}\in W^{2m-1-j-|\a|-|\b|+|\g|,\br{2m}{2m-1-|\a|},1}
\cdot L^2 \hookrightarrow L^{\br{2m}{j+|\b|-|\g|},1}\cdot L^2\hookrightarrow
L^{\br{2m}{m+1-|\g|},1}\cdot L^2\hookrightarrow L^{\br{2m}{2m-1-|\g|},1}$ and hence
\ea{l2}
  \|\<A,K\>_{j,\a,\b,\g}\|_{L^{\br{2m}{2m-1-|\g|},1}}&\le&
    c\,\|A_\a\|_{W^{2m-1-j-|\a|-|\b|+|\g|,\br{2m}{2m-1-|\a|},1}}\,\|K_\a\|_{L^2}\nn\\
  &\le&c\th\|A\|_{W^{m+1,\br{2m}{m+1},1}}
\eea
from the multiplication theorems.

Similarly, for $j=0$,
\eq{l3}
  \<A,K_0\>=\sum_{|\a|,|\b|\le m-2}\,\sum_{\g\le\b}\d^\g\<A,K\>_{0,\a,\b,\g}
\eeq
with $\<A,K\>_{0,\a,\b,\g}\in L^{\br{2m}{2m-1-|\g|},1}$ (use $\d^\a A_\a\in L^\infty$)
and the estimate
\eq{l4}
  \|\<A,K\>_{0,\a,\b,\g}\|_{L^{\br{2m}{2m-1-|\g|},1}}\le c\th\|A\|_{W^{m+1,\br{2m}{m+1},1}}.
\eeq
Now we turn our attention to $j\ge m-1$. This time, $K_j$ is a  function,
while $\d^\a D^jA_\a$ may be only a distribution. Hence we have to shift
derivatives in the opposite direction. Since $K_j\in W^{j+1-m,2}$, we can
shift at most $j+1-m$ derivatives; we will actually shift
$\min\{|\a|,j+1-m\}$. We proceed as before. In the fist case, $|\a|\le j+1-m$,
we have
\eq{l5}
  \<\d^\a D^jA_\a,K_j\>=\sum_{\g\le\a}\d^\g c_{\a\g}\<D^jA_\a,\d^{\a-\g}K_j\>.
\eeq
Since $D^jA_\a\in W^{2m-1-j,\br{2m}{2m-1-|\a|},1}\hookrightarrow
L^{\br{2m}{j-|\a|},1}$ and $\d^{\a-\g}K_j\in W^{j+1-m+|\g|-|\a|,2}\hookrightarrow
L^{\br{2m}{2m-j-1+|\a|-|\g|}}$, we have
$\<D^jA_\a,\d^{\a-\g}K_j\>\in L^{\br{2m}{2m-1-|\g|},1}$ and
\eq{l6}
  \|\<D^jA_\a,\d^{\a-\g}K_j\>\|_{L^{\br{2m}{2m-1-|\g|},1}}\le c\th\|A\|_{W^{m+1,\br{2m}{m+1},1}}.
\eeq
In the second case, $j+1-m<|\a|$, we choose any $\b\le\a$ with $|\b|=j+1-m$.
Then
\eq{l7}
  \<\d^\a D^jA_\a,K_j\>=\sum_{\g\le\b}\d^\g[c_{\b\g}
  \<\d^{\a-\b}D^jA_\a,\d^{\b-\g}K_j\>] 
\eeq
with $\d^{\a-\b}D^jA_\a\in W^{m-|\a|,\br{2m}{2m-1-|\a|},1}\hookrightarrow
L^{\br{2m}{m-1},1}$, $\d^{\b-\g}K_j\in W^{|\g|,2}\hookrightarrow
L^{\br{2m}{m-|\g|}}$, and hence $\<\d^{\a-\b}D^jA_\a,\d^{\b-\g}K_j\>\in
L^{\br{2m}{2m-1-|\g|},1}$, with the estimate \re{l6} holding also in the second
case.

Combining \re{l1}, \re{l3}, \re{l5}, and \re{l7}, we find that
\eq{l9}
  \sum_{j=0}^{2m-2}\<D^jA,K_j\>=\sum_{|\g|\le m-2}\d^\g\<A,K\>_\g\,,
\eeq
where the terms on the right-hand side denoted symbolically by $\<A,K\>_\g$
are estimated via \re{l2}, \re{l4}, and \re{l6} by
\eq{l10}
  \|\<A,K\>_\g\|_{L^{\br{2m}{2m-1-|\g|},1}}\le c\th\,\|A\|_{W^{m+1,\br{2m}{m+1},1}}.
\eeq
Note that the single summand $K_0$ in \re{hilfs} fits into the same scheme
because of
\eq{l11}
  K_0=\sum_{|\g|\le m-2}\d^\g K_0^\g\qquad\mbox{ with }K_0^\g\in L^{\br{2m}{m+1},1}
  \hookrightarrow L^{\br{2m}{2m-1-|\g|},1}.
\eeq
A similar treatment is necessary for the term $\del B\,P$ in \re{hilfs}.  
We find
\ea{l12}
  \del B\,P&=&\sum_{|\a|\le m-2}\sum_{\g\le\a}\d^\g[\del B_\a\,\d^{\a-\g}P]\nn\\
  &=&\sum_{|\g|\le m-2}\d^\g\Big(\del B_\g\,P
    +\sum_{\a>\g,|\a|\le m-2}\del B_\a\,\d^{\a-\g}P\Big)\nn\\
  &=:&\sum_{|\g|\le m-2}\d^\g[\del B_\g\,P+\<\del B,P\>_\g]
\eea
with $\del B_\a\d^{\a-\g}P\in L^{\br{2m}{2m-1-|\a|},1}\cdot W^{m-|\a|+|\g|,2}
\hookrightarrow L^{\br{2m}{2m-1-|\a|},1}\cdot L^{\br{2m}{|\a|-|\g|}}
\hookrightarrow L^{\br{2m}{2m-1-|\g|},1}$, which also implies the estimate
\ea{l13}
  \|\del B_\g\,P+\<\del B,P\>_\g\|_{L^{\br{2m}{2m-1-|\g|},1}}&\le&
    c\Big(\sum_{|\a|\le m-2}\|B_\a\|_{W^{1,\br{2m}{2m-1-|\a|},1}}\Big)\|dP\|_{W^{m-1,2}}\nn\\
  &\le&c\th\|B\|_{W^{3-m,\br{2m}{m+1},1}}.
\eea
Combining \re{l9}, \re{l11}, and \re{l12}, we rewrite \re{hilfs} as
\eq{l14}
  \sum_{|\g|\le m-2}\d^\g[\underbrace{d\Delta^{m-1}A_\g+\<A,K\>_\g+K_0^\g
    -\<\del B,P\>_\g-\del B_\g\,P}_{=:G_\g}]=0.
\eeq
Our approach to finding $A$ and $B$ is solving $[...]=0$ for every single
$\g$, now that we have more unknown functions $A_\g,B_\g$ at our disposition.
More precisely, we follow the idea from \cite{LR} to solving $[...]=0$ by
splitting it in two equations $\del[...]=0$ and $d([...]P^{-1})=0$. This gives
a system of equations
\ea{sys1}
  -\Delta^m A_\g&=&\del[\<A,K\>_\g+K_0^\g-\<\del B,P\>_\g-\del B_\g\,P],\\
\la{sys2}
  d\del B_\g&=&d[\{d\Delta^{m-1}A_\g+\<A,K\>_\g+K_0^\g-\<\del B,P\>_\g\}P^{-1}].
\eea
Note that this is a large system with two equations for every $\g$ with
$|\g|\le m-2$. We are going to choose suitable boundary conditions for each
equation separately.

For each $\g$, we define the Banach space $E_\g$ to be
$W^{2m-1,\br{2m}{2m-1-|\g|},1}$, but equipped with the norm
\[
  \|f\|_{E_\g}:=\|f\|_{W^{2m-1,\br{2m}{2m-1-|\g|},1}}
    +\th^{-1/2}\|d\Delta^{m-1}f\|_{L^{\br{2m}{2m-1-|\g|},1}}.
\]
Then we let $E:=\bigoplus_{|\g|\le m-2}E_\g$,
$E':=\bigoplus_{|\g|\le m-2}W^{2m-1,\br{2m}{2m-1-|\g|},1}$,
$L:=\bigoplus_{|\g|\le m-1}L^{\br{2m}{2m-1-|\g|},1}$,
$S:=\bigoplus_{|\g|\le m-2}W^{1,\br{2m}{2m-1-|\g|},1}$, and finally $R:=E\oplus S$.
Writing $\AA:=(A_\g)_{|\g|\le m-2}$ and $\BB:=(B_\g)_{|\g|\le m-2}$, we are looking
for $(\AA,\BB)\in R$ such that \sysa{\g} and \sysb{\g} hold for all $\g$
with $|\g|\le m-2$.

We add boundary conditions and conditions on the $dB_\g$ in order to make the
system uniquely solvable. For each $\g$, we want to solve the

{\bf Boundary value problem BVP. }{\em Solve the system made of all
\sysa{\g} and \sysb{\g} together with the Neumann boundary conditions
\eas
  \frac{\d\Delta^jA_\g}{\d\nu}&=&0\qquad\mbox{ on }\d B^{2m},\\
  \int_{B^{2m}}\Delta^jA_\g&=&0
\eeas
for all $\g$ and all $j\in\{0,\ldots,m-1\}$, the additional equations
\[
  dB_\g=0
\]
for all $\g$, as well as the boundary conditions
\[
  (B_\g)_N=0\qquad\mbox{ on }\d B^{2m}
\]
for all $\g$.}

We will solve BVP by an iteration procedure. For any $(\aa,\bb)\in R$, let
$T(\aa,\bb)=:(\aa',\bb')$ be the solution of the system
\eas
  -\Delta^ma'_\g&=&\del[\<a,K\>_\g+K_0^\g-\<\del b,P\>_\g-\del b_\g\,P],\\
  d\del b'_\g&=&d[\{d\Delta^{m-1}a_\g+\<a,K\>_\g+K_0^\g
    -\<\del b,P\>_\g\}P^{-1}],\\
  \frac{\d\Delta^ja'_\g}{\d\nu}&=&0\qquad\mbox{ on }\d B^{2m}\mbox{ for }
    j\in\{0,\ldots,m-1\},\\
  \int_{B^{2m}}\Delta^ja'_\g&=&0\qquad\mbox{ for }j\in\{0,\ldots,m-1\},\\ 
  db'_\g&=&0,\\
  (b'_\g)_N&=&0\qquad\mbox{ on }\d B^{2m},
\eeas
which must hold for all $\g$ with $|\g|\le m-2$. We further split that problem
and let $(\aa^1,\bb^1)$ be the solution of
\eas
  -\Delta^ma^1_\g&=&\del K_0^\g,\\
  d\del b^1_\g&=&d[K_0^\g P^{-1}],
\eeas
while $(\aa^2,\bb^2)$ is assumed to solve
\eas
  -\Delta^ma^2_\g&=&\del[\<a,K\>_\g-\<\del b,P\>_\g-\del b_\g\,P],\\
  d\del b^2_\g&=&d[\{d\Delta^{m-1}a_\g+\<a,K\>_\g
    -\<\del b,P\>_\g\}P^{-1}],
\eeas
both with the last four conditions similar to the system for $(\aa',\bb')$.  

By standard estimates for boundary value problems for differential forms
(see \cite{LR} for a closely related system), the systems for
$(\aa^1,\bb^1)$ and $(\aa^2,\bb^2)$ are uniquely solvable with estimates
\ea{l15}
  \|a^1_\g\|_{W^{2m-1,\br{2m}{2m-1-|\g|},1}}&\le&c\|K_0^\g\|_{L^{\br{2m}{2m-1-|\g|},1}}\,,\\
\la{l16}
  \|b^1_\g\|_{W^{1,\br{2m}{2m-1-|\g|},1}}&\le&c\|K_0^\g P^{-1}\|_{L^{\br{2m}{2m-1-|\g|},1}}\nn\\
  &\le&c\|K_0^\g\|_{L^{\br{2m}{2m-1-|\g|},1}}\,,\\
\la{l17}
  \|a^2_\g\|_{W^{2m-1,\br{2m}{2m-1-|\g|},1}}
  &\le&c\|\<a,K\>_\g-\<\del b,P\>_\g-db_\g\,P\|_{L^{\br{2m}{2m-1-|\g|},1}}\nn\\
  &\le&c\th(\|\aa\|_{E'}+\|\bb\|_S),\\
\la{l18}
  \|b^1_\g\|_{W^{1,\br{2m}{2m-1-|\g|},1}}
  &\le&c\|d\Delta^{m-1}a_\g+\<a,K\>_\g-\<\del b,P\>_\g\|_{L^{\br{2m}{2m-1-|\g|},1}}\nn\\  
  &\le&c\|d\Delta^{m-1}\aa\|_L+c\th(\|\aa\|_{E'}+\|\bb\|_S),
\eea
the latter two coming from \re{l10} and \re{l14}. The estimates show that
$T$ maps $R$ to itself.

Now we show that $T$ is a contraction on $R$. Let
$(\aa,\bb), (\aas,\bbs)\in R$. Since $(\aa^1,\bb^1)=(\aas^1,\bbs^1)$, we have
estimates similar to \re{l17} and \re{l18},
\eas
  \|\aa'-\aas'\|_{E'}&=&\|\aa^2-\aas^2\|_{E'}\\
  &\le&c\th(\|\aa-\aas\|_{E'}+\|\bb-\bbs\|_S),\\
  \|\bb'-\bbs'\|_S&=&\|\bb^2-\bbs^2\|_S\\
  &\le&c\|d\Delta^{m-1}(\aa-\aas)\|_L+c\th(\|\aa-\aas\|_{E'}+\|\bb-\bbs\|_S)\\
  &\le&c\th^{1/2}\|\aa-\aas\|_E+c\th(\|\aa-\aas\|_{E'}+\|\bb-\bbs\|_S),
\eeas
from which we infer
\[
  \|(\aa',\bb')-(\aas',\bbs')\|_R\le c\th^{1/2}\|(\aa,\bb)-(\aas,\bbs)\|_R.
\]
Hence $T$ is a contraction provided that $\th$ has been chosen small enough.  
It therefore has a fixed point $(\AA,\BB)\in R$. The mapping $T$ has been set
up in such a way that this fixed point solves $BVP$.

By estimates corresponding to
\re{l15}--\re{l18} and $\AA'=\AA$, $\BB'=\BB$, we have
\eas
  \|(\AA,\BB)\|_R&=&\|\AA^1\|_E+\|\BB^1\|_S+\|\AA^2\|_E+\|\BB^2\|_S\\
  &\le&c\th^{-1/2}\|K_0\|_L+c\|K_0\|_L+c\th^{1/2}(\|\AA\|_E+\|\BB\|_S)\\
  &\le&c\th^{1/2}(1+\|(\AA,\BB)\|_R).
\eeas
Again if $\th>0$ has been chosen small enough, then this implies
$\|(\AA,\BB)\|_R\le c\th^{1/2}$. Letting $A:=\sum_{|\a|\le m-2}\d^\a A_\a$ and
$B:=\sum_{|\a|\le m-2}\d^\a B_\a$, we find that the norms in \re{ABklein} are
easily seen to be controlled by $\|(\AA,\BB)\|_R$, thus \re{ABklein} is proven.

Now that we have solved BVP, we need to show that \re{hilfs} is solved by
our $A$ and $B$, and we have already seen that this holds once we have
proven that \re{sys1} and \re{sys2} (which are $\del G_\g=0$ and
$d(G_\g P^{-1})=0$) together with the boundary data from $BVP$
imply $G_\g=0$, with $G_\g$ defined in \re{l14}. 
Close to the boundary, the vanishing of $K_j$ and $dP$ simplifies the
expression for $G_\g$, giving
\[
  G_\g=d\Delta^{m-1}A_\g-\del B_\g
\]
there. This means that the normal component of $G_\g$ vanishes, because we have
Neumann boundary data for $\Delta^{m-1}A$, and $(B_\g)_N=0$ on $\d B^{2m}$
also implies $(\del B_\g)_N=0$.

We have $G_\g\in L^{\br{2m}{2m-1-|\g|},1}(\R^{n\times n}\otimes\wedge^1\R^{2m})$.
A standard version of the Hodge decomposition (generalized to Lorentz
spaces) says that $\del G_\g=0$ on $B^{2m}$ and and $(G_\g)_N=0$ on
$\d B^{2m}$ imply $G_\g=\del C_\g$ for some
$C_\g\in W^{1,\br{2m}{2m-1-|\g|},1}(\R^{n\times n}\otimes\wedge^2\R^{2m})$ with
$dC_\g=0$ on $B^{2m}$ and $(C_\g)_N=0$ on $\d B^{2m}$. We combine the latter
two equations with $d(G_\g P^{-1})=0$ to the system
\[
  d((\del C_\g)P^{-1})=0,\qquad dC_\g=0,\qquad
  ((C_\g)_N)_{|\d B^{2m}}=0.
\]
By Lemma \ref{TODO}, this implies $C_\g=0$ once $\th$ has been chosen small
enough to make $\|d(P^{-1})\|=\|dP\|$ sufficiently small. Then also
$G_\g=0$, which implies \re{hilfs} and proves the lemma.\qed

\section{Regularity}

From Theorem \ref{th} (iii), we now easily infer that weak solutions of
our differential equation \re{eq} are continuous.

\begin{theorem}\la{co}
  Assume $U$ is a bounded open set in $\R^{2m}$ and $u\in W^{m,2}(U,\R^n)$
  is a weak solution of \re{eq} on $U$ with coefficient functions
  in the same spaces as in Theorem \ref{th} (but without smallness
  condition). Then $u$ is continuous on $U$.
\end{theorem}

{\bf Proof.} Fix $x_0\in U$. For any $r>0$ small enough that
$B_r(x_0)\subseteq U$ define $u_r:B^{2m}\to\R^n$ by
$u_r(x):=u(x_0+rx)$. Then $u_r$ solves an equation on $B^{2m}$ of the form
\re{eq}, but with $V_k,w_k$ replaced by $V_{k,r}(x):=r^{2m-2k-2}V_k(x_0+rx)$
and $w_{k,r}(x):=r^{2m-2k-2}w_k(x_0+rx)$. We have
$\|V_{k,r}\|_{W^{2k+1-m,2}(B^{2m})}\le\|V_k\|_{W^{2k+1-m-2}(B_r(x_0))}$ and
$\|w_{k,r}\|_{W^{2k+2-m,2}(B^{2m})}\le\|w_k\|_{W^{2k+2-m-2}(B_r(x_0))}$, both of which
can be made arbitrarily small if $r>0$ is chosen small. Similar estimates
hold for $\eta$ and $F$, and hence the smallness condition $\th<\th_0$
from Theorem \ref{th} is fulfilled for the differential equation for $u_r$
once $r$ has been chosen sufficiently small. But then $u_r$ is continuous
on $B_{1/16}$, which means $u$ is continuous on $B_{r/16}(x_0)$.\qed

Of course, in many special cases, we can expect much more than just continuity
of weak solutions. As was proven in \cite[Section 7]{GS}, extrinsically or
intrinsically polyharmonic maps are even smooth in $2m$ dimensions, and,
for a large class of equations, H\"older continuity of weak solutions
implies smoothness. But here, we consider a very general equation
with rather irregular coefficients, so maybe we cannot expect much
regularity in general.

\vfill

\begin{center}
\scriptsize
Fakult\"at f\"ur Mathematik, Universit\"at Duisburg-Essen, D-47048 Duisburg,
Germany\\
{\tt frederic.de-longueville@uni-due.de}\\
\strut\\
Fakult\"at f\"ur Mathematik, Universit\"at Duisburg-Essen, D-45117 Essen, 
Germany.\\
{\tt andreas.gastel@uni-due.de}
\end{center}

\end{document}